\documentclass[reqno,12pt]{article}
\usepackage{amsmath, amsthm, amscd,amsfonts, amssymb,hyperref,floatrow}
\usepackage{graphicx, color,dsfont,xcolor}

\usepackage{authblk}

\numberwithin{equation}{section}

\newtheorem{theorem}{Theorem}[section]
\newtheorem{lemma}[theorem]{Lemma}
\newtheorem{proposition}[theorem]{Proposition}

\newtheorem{rem}[theorem]{Remark}

\newtheorem{hypothesis}{Hypothesis}

\DeclareMathSymbol{\leqslant}{\mathalpha}{AMSa}{"36} 
\DeclareMathSymbol{\geqslant}{\mathalpha}{AMSa}{"3E} 
\DeclareMathSymbol{\eset}{\mathalpha}{AMSb}{"3F}     
\renewcommand{\leq}{\;\leqslant\;}                   
\renewcommand{\geq}{\;\geqslant\;}                   

\def \C{ \mathbb  C }
\def \H{ \mathbb  H }
\newcommand{\R}{\mathbb{R}}

\newcommand{\ind}{\mathds{1}}

\def \E{ \mathbb E  }

\newcommand \be  {\begin{equation*}}
\newcommand \bea {\begin{eqnarray} \nonumber }
\newcommand \ee  {\end{equation*}}

\newcommand \ba  {\begin{align}}
\newcommand \ea  {\end{align}}

\definecolor{remi}{rgb}{0,0,0}

\begin{document}
\title{Eigenvector dynamics under free addition}
\author[1,2]{Romain Allez}
\author[2]{Jean-Philippe Bouchaud}

\affil[1]{Weierstrass Institute\\ Mohrenstr. 39, 10117 Berlin, Germany.}
\affil[2]{Capital~Fund~Management. 
23-25, rue de l'Universit\'e, 75\,007 Paris. France.
}

\maketitle 

\begin{abstract}
We investigate the evolution of a given eigenvector of a symmetric (deterministic or random) matrix 
under the addition of a matrix in the Gaussian orthogonal ensemble. 
We quantify the overlap between this single vector with the eigenvectors 
of the initial matrix and identify precisely a ``Cauchy-flight'' regime. 
In particular, 
we compute the local density of this vector in the eigenvalues space of the initial matrix. 
Our results are obtained in a non perturbative setting and are derived using the ideas of 
[O. Ledoit and S. P\'ech\'e, Prob. Th. Rel. Fields, {\bf 151} 233 (2011)]. 
Finally, we give a robust 
derivation of a result obtained in [R. Allez and J.-P. Bouchaud,  Phys. Rev. E {\bf 86}, 046202 (2012)]
to study eigenspace dynamics in a semi-perturbative regime.
\end{abstract}

\section{Introduction}

The dynamics of eigenvalues induced by the addition of 
free random matrices in the Gaussian
orthogonal ensemble has been first studied by Dyson in his 1962 paper \cite{dyson}. 
The movement of the eigenvalues is characterized in terms of a stochastic differential system, the so called 
{\it Dyson Brownian motion}. The eigenvalues evolve as particles of a Coulomb gas with 
electrostatic repulsion, confined in a quadratic potential and subject to a thermal noise.
In the limit of large dimensions (matrix sizes), the evolution of the spectral density has also been studied 
in \cite{rogers-shi} (see also \cite{biane,florent} and \cite{cepa,jp-alice} for related models).      

For the eigenvectors, their evolution in finite dimension is also given by a stochastic differential system
which depend on the non colliding trajectories of the eigenvalues (see \cite{agz}). 
In this paper, we are interested in quantifying the evolution of the eigenvectors in the limit of large dimension. 
Our approach uses the idea of \cite{sandrine} who introduced a very interesting quantity (see Eq. \eqref{Theta_N.def} below) 
for the study of eigenvectors. This enables us to compute the local density of a given state (eigenvector) 
of the matrix after the addition of the free Gaussian matrix, in the eigenvalues space of the initial matrix. 

The paper is organized as follows. In section \ref{def}, we define the model and give the main notations. 
 In section \ref{finiteN}, we first consider the evolution of the eigenvectors induced by the addition of a small Gaussian matrix 
 when the dimension $N$ of the matrices are small.   
Our main Theorem \ref{main} appears in section \ref{infiniteN} and is concerned with 
the convergence of the quantity \eqref{Theta_N.def} introduced in \cite{sandrine}. The proof of Theorem \ref{main}
can be found in section \ref{proof}. 
In section \ref{semi-circular-initial}, we specialize our results in a natural case where the computations are explicit
and we find that in a particular regime, the eigenvalue dynamics can be precisely described as a ``Cauchy flight''.
We also check numerically our results in the case of a initial random matrix in the Gaussian orthogonal ensemble.  
We then revisit in section \ref{stability} the main result of \cite{vectors} 
on the dynamics of eigenspace under free addition, and prove that it is indeed exact beyond the perturbative regime.

{\bf Acknowledgments}   
We are grateful to Yan Fyodorov, Denis Ullmo and Matthieu Wyart who asked questions about \cite{vectors} which led to the work presented here. 
We would like to thank also Sandrine P\'ech\'e for interesting discussions.  
We thank the anonymous referee for his careful reading of previous versions of this paper and for his valuable suggestions.  

RA received funding from the European Research Council under the European
Union's Seventh Framework Programme (FP7/2007-2013) / ERC grant agreement nr. 258237 and thanks the Statslab in DPMMS, 
Cambridge for its hospitality.

\section{Definition and main notations}\label{def}
Let $A$ be a symmetric deterministic $N \times N$ matrix. By the spectral Theorem, $A$ is diagonalizable in an orthonormal basis of $\R^N$. 
We suppose that the eigenvalues $a_1,\dots, a_N$ of $A$ are all distinct and indexed in increasing order as
\begin{align}\label{assumptionA}
a_1 < a_2 < \cdots < a_N\,. 
\end{align}

Let $(H(t))_{t\geq 0}$ be a symmetric Brownian motion, i.e. a symmetric diffusive matrix process constructed from a family of independent real
Brownian motions $B_{ij}(t)$, $1 \leq i \leq j \leq N$ as follows 
\begin{equation*}
H_{ij}(t)=\begin{cases} \sqrt{\frac{1}{N}} \, B_{ij}(t) & \quad {\rm if} \, i < j \,, \\ \sqrt{\frac{2}{N}} \,  B_{ii}(t) & \quad {\rm if} \, i=j \,.  \end{cases}
\end{equation*}
The process $H(t)$ is rotationally invariant at all time $t\geq 0$, 
in the sense that for all $O$ in the orthogonal group $\mathcal{O}_N$, 
the conjugate matrix $O \,H(t)\, O^\dagger$ has the same law as the matrix $H(t)$. 

 Now we define a Hermitian matrix process $(M(t))_{t\geq 0}$ by setting  
\begin{equation}\label{eq.Mt}
M(t) = A + H(t) \,. 
\end{equation}
The matrix $M(t)$ may be regarded as a noisy perturbation of the matrix $A$, which encodes the {\it true} information. 
The eigenvalues of $M(t)$ will be denoted in increasing order as 
\begin{align*}
\lambda_1(t) \leq \cdots \leq \lambda_N(t)\,.
\end{align*}
The aim of this paper is to quantify the relationship between the eigenvectors of the matrix $M(t)$ with the ones 
of the initial matrix $M(0)=A$. 
In particular, we consider one given eigenvector of the matrix $M(t)$ denoted as $\psi_i^t$ and we want to 
compute, in the limit of large dimension $N$, 
its projections  
on the (orthonormal) eigenvectors of $A$ denoted in the following as 
$\phi_1, \cdots, \phi_N$. The scalar products between the vectors $\psi_i^t$ and $\phi_j$ 
are also called {\it overlaps} between $\psi_i^t$ and
$\phi_j$ and denoted $\langle \psi_i^t |\phi_j\rangle$. 
Because the matrix $H(t)$ is rotationally invariant for all $t$, we can (and will) suppose 
with no loss of generality that the matrix $A$ is {\it diagonal}.

In order to study the limit $N\to +\infty$, we 
need to take a few (natural) assumptions on the spectrum of the matrix $A:=A_N$. 
 \begin{hypothesis}\label{hypothesis}
We suppose that the empirical spectral density of the matrix $A$, defined as
\begin{align*}
\mu_0^N:=\frac{1}{N} \sum_{i=1}^N \delta_{a_i}
\end{align*}
where $\delta_x$ is the Dirac measure in $x$,
converges in the space of probability measures (equipped with the topology of weak convergence) as 
$N\to \infty$ to $\mu_0(dx):=\rho_0(x)\, dx$, where $\rho_0:\R\to \R_+$ is a continuous function. 
In fact, we ask for a little more, making the assumption that the $a_i$ are allocated smoothly 
on the quantiles of the probability density $\rho_0(x)$, i.e.
according to $a_i=a(\frac{i}{N+1}),i=1,\dots,N$ where
the function $a:]0,1[\to \R$ is continuous, strictly increasing and such that for all $x\in ]0,1[$,
\begin{align*}
\int_{-\infty}^{a(x)} \rho_0(y) \, dy = x\,. 
\end{align*}
This definition implies that the push-forward measure $a^{-1}(\mu_0)$ (respectively 
$a^{-1}(\mu_0^N)$) of $\mu_0$ (resp. $\mu_0^N$) 
by the function $a^{-1}: \R \to (0,1)$  is 
the uniform measure on $[0,1]$ (resp. the discretized uniform measure $\frac{1}{N} \sum_{i=1}^N \delta_{i/(N+1)}$).  
\end{hypothesis}
Note that if $\rho_0$ has unbounded support, then $a(x)\downarrow -\infty$ as $x\downarrow 0$ and $a(x)\uparrow +\infty$ as $x\uparrow 1$. 
If $\rho_0$ has a compact support, then the function $a$ is bounded on $]0;1[$ and converges when $x\downarrow 0$ and $x\uparrow 1$.

\section{Dyson Brownian motion and eigenvectors dynamics in finite dimension $N$} \label{finiteN}

We first discuss asymptotic formulas for the scalar products between the eigenvectors of $M(t)$ 
with those of the matrix $M(0)=A$ in the perturbative limit $t\to 0$ when the dimension $N$ is {\it fixed}.
In the next section, we shall investigate the limit $N\to \infty$, with possible scaling relations between $t$ and $N$.    

In the present case, the idea is to write the evolution equations for the eigenvalues and eigenvectors 
of the Hermitian matrix $M(t)$.
It was first established by Dyson in \cite{dyson} that the eigenvalues follow the stochastic differential system 
\begin{align}\label{sde.ev}
d\lambda_i(t) = \sqrt{\frac{2}{N}}\, dB_{i}(t) + \frac{1}{N} \sum_{j \neq i } \frac{dt}{\lambda_i-\lambda_j}\,,
\end{align} 
where the $B_i$ are independent Brownian motions and  
with the initial conditions $\lambda_i(0)=a_i, i =1, \cdots,N$. 
Let us simply mention that the electrostatic repulsion (last term of \eqref{sde.ev}) is strong enough 
to prevent any collision between the eigenvalues so that the stochastic differential system has a well defined 
and continuous solution in the It\^o's sense \cite{agz}. Recalling the assumption \eqref{assumptionA} on the
location of the $a_i=\lambda_i(0)$ at the initial time, we can conclude that for all $t\geq 0$, 
\begin{align*}
\lambda_1(t) < \lambda_2(t) < \dots < \lambda_N(t)\,. 
\end{align*} 
Towards a physical picture, we can see the process $(\lambda_1,\dots,\lambda_N)$ 
as a one dimensional repulsive Coulomb gas of $N$ positively charged particles, all subject to a thermal noise $dB_i(t)$. 

For the evolution of the eigenvectors $(\psi_1^t,\dots, \psi_N^t)$ as a function of time $t$, 
the situation is slightly more tricky because the $\psi_i^t$
are all determined up to a sign $\pm 1$. Nevertheless, 
we can prove following \cite[Proof of Theorem 4.3.2]{agz} (see also \cite{wilkinson,alice})
that there exists a continuous (with respect to time) version of the process $(\psi_1,\dots,\psi_N)$ 
which evolves according to 
\begin{align}\label{sde-vectors}
d\psi_i^t = - \frac{1}{2N} \sum_{j \neq i} \frac{dt}{(\lambda_i-\lambda_j)^2}\,  \psi_i^t
+ \frac{1}{\sqrt{N}}\sum_{j \neq i} \frac{dW_{ij}(t)}{\lambda_i-\lambda_j} \, \psi_j^t\,,
\end{align}
with the initial conditions $\psi_i^0=\phi_i$ for $i=1,\dots,N$, and 
where the family of real Brownian motions $W_{ij}, 1 \le i , j \le N$ is such that  
\begin{itemize}
\item the $W_{ij}, i \le j$ are mutually independent;
\item the $W_{ij}$ for $i > j$ are defined by symmetry $W_{ij} = W_{ji}$;
\item the $W_{ij}, i \le j$ are independent of the Brownian motions $B_i$
driving the stochastic differential system of the eigenvalues \eqref{sde.ev}.
\end{itemize}
The independence between the $W_{ij}, i \le j$ and the $B_i$ allows us to freeze the trajectories of the eigenvalues  
and then to study the eigenvectors dynamics with this realization of the eigenvalues path. 

With this description of the evolution of the eigenvalues and eigenvectors processes,
we easily deduce the following proposition.  
\begin{proposition}\label{proposition-finite-N}
Let $N$ be a fixed integer and $a_1 < a_2 <\cdots < a_N$ be the eigenvalues of the diagonal matrix $A$. 
Then, the (mean) overlaps between the eigenvectors of $M(t)$ defined in \eqref{eq.Mt} with those 
of $M(0)=A$ satisfy the following asymptotic expansions when $t\to 0$, 
\begin{align}\label{overlap.diag}
\E\left[\langle\psi_i^t | \phi_i \rangle \right] = 1 - \frac{t}{2N} \sum_{j \neq i } \frac{1}{(a_i-a_j)^2} + o(t)\,,
\end{align} 
and for $i\neq j$,
\begin{align}\label{overlap.nondiag}
\E\left[\langle\psi_i^t | \phi_j \rangle^2\right] = \frac{t}{N} \frac{1}{(a_i-a_j)^2} +o(t)  \,.
\end{align}
\end{proposition}

{\it Proof.}
We fix $T>0$ and 
denote by $\E_\lambda$ the expectation conditionally on the $\sigma$ algebra generated by 
the eigenvalues trajectories $(\lambda_i(t), i \in \{1,\dots,N\} ,0\le t\le T)$. 
From \eqref{sde-vectors}, we see using the independence between the $(W_{ij},i \le j)$  and the $\lambda_i$ that 
\begin{align}\label{diff-eq-scalar-product}
\E_\lambda[\langle\phi_i|\psi_i^t \rangle] = 1 -\frac{1}{2N} \sum_{j \neq i} \int_0^t 
\frac{\E_\lambda[\langle\phi_i|\psi_i^s\rangle]}{(\lambda_i-\lambda_j)^2} ds\,. 
\end{align} 
We can solve \eqref{diff-eq-scalar-product} explicitly and obtain (using $\psi_i^0=\phi_i$ and $\langle \phi_i|\phi_i\rangle=1$) 
\begin{align*}
\E_\lambda[\langle\phi_i|\psi_i^t \rangle] = \exp\left(-\frac{1}{2N} \sum_{j\neq i} \int_0^t \frac{ds}{(\lambda_i(s)-\lambda_j(s))^2}\right)\,. 
\end{align*}
The asymptotic expansion \eqref{overlap.diag} follows from the fact that $\lambda_i(t)\to a_i$ almost surely as $t\to 0$.  
 
Towards \eqref{overlap.nondiag}, we use again \eqref{sde-vectors} and noting that 
$\langle \psi_i^0 |\phi_j\rangle = \langle \phi_i|\phi_j\rangle$ 
for $i\neq j$, we obtain 
\begin{align}\label{diff-eq-off-diag}
\E_\lambda[\langle  \psi_i ^t|\phi_j\rangle^2] &=\frac{1}{N}\int_0^t \sum_{k \neq i} \frac{\E_\lambda[\langle \psi_k^s|\phi_j \rangle^2  ] - \E_\lambda[\langle \psi_i^s|\phi_j \rangle^2  ] }{(\lambda_k-\lambda_i)^2} ds\,.
\end{align}
The asymptotic expansion \eqref{overlap.nondiag} now follows from the almost sure convergences as $s\to 0$ 
of  $\langle \psi_k^s |\phi_j \rangle$ towards $1$ if $k= j $ and $0$ if $k \neq j$ and of $\lambda_k(s)$ to $a_k$.  
\qed

\section{Eigenvectors dynamics in the large $N$ limit}\label{infiniteN}
We now consider the eigenvector dynamics problem of quantifying the relation between the eigenvectors of the 
Hermitian matrix $M(t)$ with those of $A$ but in the limit of large dimension $N\to \infty$. 
In this context, we shall distinguish three regimes   
characterized by different scaling relations between the two parameters $t$ and $N$:
\begin{itemize}
\item The first regime is the {\it perturbative regime} where $t:=t_N$ is scaling with the dimension $N$ and satisfies
\begin{align*}
N t_N\longrightarrow 0\,, 
\end{align*}
when $N\to +\infty$. 
In this regime, $t_N$ is in fact much smaller than the typical spacing between two consecutive eigenvalues (level spacing), which is of order $1/N$ 
in our setting. This allows {\it perturbation theory}  
to be applied in the eigenvalues problem of the matrix $M(t)$, seen as a perturbation 
of the matrix $A$, and to derive approximations identical  to \eqref{overlap.diag} and \eqref{overlap.nondiag}. 
This regime has been studied in great details in random matrix theory and in the context of quantum mechanics 
\cite{fernandez}.  
 \item The second regime is {\it semi-perturbative} and concerns values of $t:=t_N$ 
again scaling with the dimension $N$ such that
\begin{align*}
t_N\longrightarrow 0\,,
\end{align*}
but which are not necessarily small compared to the levels spacing of order $1/N$. 
This regime concerns many applications (see e.g. \cite{sandrine,vectors} in 
the context of covariance matrices and applications to finance) and basic perturbation theory 
does not permit one to rigorously extend the validity of Eq. \eqref{overlap.diag} and \eqref{overlap.nondiag} 
to this regime. Our main result Theorem \ref{main} permits us to do so (see the discussion after Theorem \ref{main}).   
 \item The third regime is {\it non perturbative}: $t$ is fixed independent of the dimension $N$ going to $+\infty$.  
\end{itemize}

The question we ask is: How to modify formulas \eqref{overlap.diag} and \eqref{overlap.nondiag} 
 in the second and third regimes in the large $N$ limit ? 
Because the family $\{\phi_j,1\leq j \leq N\}$ forms an orthonormal basis of $\R^N$, we have the normalization constraint
\begin{equation*} 
\sum_{j=1}^N \langle \psi_i^t |\phi_j\rangle^2 = 1\,,
\end{equation*} 
and shall therefore investigate the convergence of the renormalized overlaps 
$N\, \E[\langle \psi_i^t | \phi_j\rangle^2]$ for $i\neq j$ as $N\to \infty$. Those scalar products are in fact related to the mean
{\it local density of the state} $|\psi_i^t\rangle$ which is defined as the probability measure $\nu_i$, supported 
on the eigenvalues of $A$, 
\begin{equation*}
\nu_i(da) = \sum_{j=1}^N \E[\langle \psi_i^t |\phi_j\rangle^2] \,\, \delta_{a_j}(da)\,.
\end{equation*} 
 In other words, the aim of this paper is to compute the 
local density $\nu_i$ of the eigenvector $|\psi_i^t\rangle$ in the large $N$ limit.   

The interesting quantity for our purpose is the bivariate cumulative distribution function $\Phi$ associated 
to the weights $N \, \E[\langle \psi_i^t|\phi_j\rangle^2]$ defined for $\lambda,\alpha \in \R$ by  
\begin{align}\label{bivariate-distrib}
\Phi_N(\lambda, \alpha) = \frac{1}{N} \sum_{i,j=1}^{N} \E[ \langle \psi_i^t|\phi_j\rangle^2]  \, {\ind}_{\{\lambda_i^t \leq \lambda\}} \,  {\ind}_{\{a_j \leq \alpha \}}\,.  
\end{align}
Note that this function $\Phi$ has indeed the properties of a bivariate cumulative distribution function since  
\begin{itemize}
\item it is right continuous with left-hand limits;
\item it is nondecreasing in each of its argument;
\item it satisfies $\lim\limits_{\lambda\rightarrow-\infty, a \rightarrow -\infty} \Phi(\lambda, a)= 0 $ and 
$\lim\limits_{\lambda\rightarrow+\infty, a \rightarrow +\infty}\Phi(\lambda, a) = 1$. 
\end{itemize}

Before presenting our results on the convergence as $N\to \infty$ of the bivariate cumulative distribution $\Phi_N(\lambda,\alpha)$ 
which shall directly lead us to asymptotic estimates for the overlaps $N \E[\langle \psi_i^t| \phi_j\rangle^2 ]$ for $i\neq j$ and for the local density of states, we state a Theorem due to Shlyakhtenko \cite{shlyakhtenko} 
(see also \cite{pastur,zee} for similar results)
on the convergence of the empirical eigenvalue distribution of the matrix $M_t$ defined in \eqref{eq.Mt}. 
Let us recall the definition of the Stieltjes transform $G_\mu(z)$ of a probability measure $\mu$ on $\R$ defined 
on the upper half plane $\H :=\{z \in \C: \Im z >0\}$ by 
\begin{align*}
G_\mu(z)=\int_\R \frac{\mu(dx)}{x-z}\,. 
\end{align*} 
The Stieltjes transform is frequently used in random matrix theory for the study of empirical spectral densities in the large $N$ limit. 
A measure $\mu$ is characterized by its Stieltjes transform, which is an analytic function $G: \H \to \H$. 
We have the following inversion formula valid for any measure $\mu$ on $\R$ and $x<y$,
\begin{align}\label{inversion-formula}
\lim_{\varepsilon \downarrow 0}\int_x^y \Im \, G_\mu(\lambda+i \, \varepsilon) \, d\lambda
= \pi\, \mu(x;y) + \frac{\pi}{2} \left(\mu(\{y\}) - \mu(\{x\})\right)\,,  
\end{align}
where $\Im z$ denotes the imaginary part of $z\in \C$. 
If $\mu$ is a probability measure, its Stieltjes transform $G_\mu(z)$ behaves as $-1/z$ when $|z|$ goes to infinity.


\begin{proposition}[Shlyakhtenko, \cite{shlyakhtenko}]\label{Shlyakhtenko.th}
For $z\in \H$, set $R_t^N(z)= (M_t - z I)^{-1}$ and introduce the complex measure 
defined on the unit interval $[0,1]$ as 
\begin{equation}\label{def.sigma}
\sigma_t^N(z,dx)=\frac{1}{N} \sum_{i=1}^N \E[R_t^N(z)_{ii}] \,\, \delta_{\frac{i}{N}}(dx)\,. 
\end{equation}
Then, the complex measure $\sigma_t^N(z,dx)$ converges weakly to a complex measure with density $G_t(z,x) dx$, 
where $G_t(\cdot,x):\H \to \H$ is the unique analytic function 
such that, for all $z\in \H$ and $x\in [0,1]$,
\begin{equation}\label{eq.pt-fixe}
G_t(z,x) =  \frac{1}{a(x)-z - t \int_{0}^1 G_t(z,y) dy} \,. 
\end{equation}
\end{proposition}

\begin{rem}
The Stieltjes transform $G_{\mu_t^N } (z)$ of the mean empirical spectral measure $\mu_t^N$ of $M_t$, defined as 
\begin{align*}
\mu_t^N(d\lambda) = \frac{1}{N} \sum_{i=1}^N \E[\delta_{\lambda_i}(d\lambda)]\,,
\end{align*}
can be recovered from 
$\sigma_t^N(z,dx)$ by the formula 
\begin{align*}
G_{\mu_t^N } (z) = \sigma_t^N(z,[0,1])\,. 
\end{align*} 
From this observation, it is easy to see that $\mu_t^N$ converges weakly  
when $N\rightarrow \infty$ to the probability 
measure $\mu_t$ associated to the Stieltjes transform 
\begin{equation*}
G_{\mu_t}(z) = \int_0^1 G_t(z,x) \, dx\,. 
\end{equation*} 
\end{rem}

Theorem \ref{Shlyakhtenko.th} can be proved using the Schur complement formula, which permits one to obtain the 
equation \eqref{eq.pt-fixe} satisfied by the limit points along subsequences of the complex measure $\sigma_t^N$,
and then the fixed point theorem to show the uniqueness of the analytic function $G_t(\cdot,x)$ satisfying \eqref{eq.pt-fixe} 
(see \cite{gerard-alice,remi-vincent} where such a route was used). 
In \cite{shlyakhtenko}, Shlyakhtenko proves in fact a more general result, covering the case 
of band random matrices, with a proof using the theory of free probability.  

Finally, using the work of Biane \cite{biane} {\it On the free convolution with a Semi-circular distribution}, we know that the limiting spectral distribution $\mu_t(d\lambda)$ admits a smooth 
density $\rho_t(\lambda)$ with respect to Lebesgue 
measure. In \cite[Corollary 2]{biane}, Biane gives an analytic formula for the density $\rho_t$ as a function of the 
Stieltjes transform $G_{\mu_0}$ of the initial spectral density $\mu_0$, proving 
the convergence of the Stieltjes transform near the real axis:
\begin{equation*}
\lim\limits_{\eta\rightarrow 0_+} G_{\mu_t}(\lambda+i \eta) = H_{\rho_t}(\lambda)+ i \pi \rho_t(\lambda)
\end{equation*}
where $H_{\rho_t}$ is the Hilbert transform of the probability density $\rho_t$.  
The work of Biane is motivated by the theory of
 free probability, which was  
originally introduced by Voiculescu in \cite{voiculescu1} (see also \cite{vdn,hiai})  
as a new theory of probability for non-commuting random variables. In this context, {\it freeness}
or {\it free independence} is the analogue of the classical notion of independence. 
Using the connections between random matrices and free probability later established in \cite{voiculescu2},  
the limiting spectral measure $\mu_t$ may be seen as the {\it free convolution} between the initial probability measure 
$\mu_0$ and the semi-circular distribution of variance $t$, $\lambda_t(dx):=\frac{1}{2\pi t} \sqrt{4t^2- x^2}$. This operation
between two real measures is usually denoted $\boxplus$. In the present case, we have 
\begin{align*}
\mu_t = \mu_0 \boxplus \lambda_t\,. 
\end{align*}
This connection with random matrices motivates the title of our article: while the addition of
the Gaussian matrix $H(t)$ induces (in the large $N$ limit) 
a free convolution of the initial spectral measure $\mu_0$ by 
 the semi-circular distribution of variance $t$, we study the relation between the eigenvectors 
 of the matrix $M(t)=A+H(t)$ with those of the matrix $M(0)=A$ at the initial time.


We are now ready to state our main result on the convergence of the bivariate cumulative distribution 
$\Phi_N(\lambda,\alpha)$. 

\begin{theorem} \label{main}
Let $t>0$ and $A:=A_N$ a $N\times N$ symmetric matrix such that hypothesis \ref{hypothesis} holds
for some general initial probability density $\rho_0$. 
We consider the random matrix $M(t)=A+ H(t)$ of the matrix $A$, where $(H(t))_{t\geq 0}$ is a symmetric Brownian 
motion.  

Then, the bivariate cumulative distribution $\Phi_N$  defined in \eqref{bivariate-distrib}
converges point wise as $N\to +\infty$ to a bivariate cumulative distribution $\Phi$ given by 
\begin{align*}
\Phi(\lambda,\alpha) =
\int_{-\infty}^\lambda d\xi \, \rho_t(\xi) \int_{-\infty}^\alpha dx\,  \rho_0(x)\,  \frac{t}{(x-\xi-t \, H_{\rho_t}(\xi))^2 + t^2 \pi^2 \rho_t(\xi)^2}\,. 
\end{align*}
\end{theorem}  

As we will see, Theorem \ref{main} permits to compute the asymptotic overlaps 
$ \E[\langle \psi_i^t|\phi_j\rangle^2]$ of any eigenvector $|\psi_i^t\rangle$ of the matrix $M_t$ at time $t$ 
with the eigenvectors of the initial matrix $M_0=A$ for any time $t$. 
A previous heuristic attempt to compute those overlaps appeared in \cite{wilkinson}, with a different result (see below).  
Theorem \ref{main} also permits one to compute the mean local density $\nu_i$ in the $A$-eigenvalues space
of the state $|\psi_i^t\rangle$ at time $t$.  In addition we will explain how it enables us to extend the domain of validity 
of our former results on the eigenspace dynamics under free addition obtained in \cite{vectors}.
 
Theorem \ref{main} can be seen as the counterpart of \cite[Theorem 3]{sandrine} which quantifies
 the relationship between the eigenvectors of 
the population covariance matrix with those of the empirical (or sample) 
covariance matrix. 

One remarkable feature of Theorem \ref{main} is that it quantifies the relationship 
between the eigenvectors of the initial matrix $M_0=A$ 
and the eigenvectors of the matrix $M_t$, even in the {\it non perturbative} third regime described above 
(where $t>0$ is independent of $N$).   

Our proof of Theorem of \ref{main} uses the ideas of \cite{sandrine} to quantify 
the relationship between sample and population eigenvectors. 

Theorem \ref{main} enables to compute the overlaps of the vector $|\psi_i^t\rangle$  
with the initial ($t=0$) eigenvectors $|\phi_j\rangle$ for $j\neq i$ of the matrix $M_0=A$. 
Roughly speaking, when $N\rightarrow \infty$, we have
\begin{align}
 \E[\langle \psi_i^t | \phi_j \rangle^2 ] = \frac{1}{N}
\frac{t}{(a_j-\lambda_i^t -t H_{\rho_t}(\lambda_i^t))^2 + t^2 \pi^2 \rho_t(\lambda_i^t)^2} + o(\frac{1}{N})\,. \label{overlap}
\end{align} 
As mentioned before, formula \eqref{overlap} is valid for $t>0$ independent of $N$, which is very large (third regime).

But when $t\to 0$, $\lambda_i(t)\to a_i$ almost surely and
if $t$ is now itself very small,  
formula \eqref{overlap} can be simplified, for any pair of indices $(i,j)$ such that the eigenvalues $a_i,a_j$ remain 
separated by a macroscopic spacing (i.e. such that $|a_j-a_i|=O(1)$ 
does not vanish for large $N$),   as
 \begin{align}\label{perturb.2.n}
 \E[\langle \psi_i^t | \phi_j \rangle^2 ] = \frac{t}{N}
\frac{1}{(a_j-a_i)^2} + o(\frac{1}{N})\,.  
\end{align}   
For such pair $(i,j)$, \eqref{perturb.2.n} extends the perturbation equations \ref{overlap.nondiag} of 
Proposition \ref{proposition-finite-N}
(which was valid only for $t\to 0$ and $N$ fixed) 
to values of $t$ much smaller than $1$ but not necessarily negligible compared to $1/N$, 
for example such that $t:=t_N=1/N^{\alpha}, \alpha \in (0;1]$,  which correspond to the semi-perturbative regime. 

One can also use Theorem \ref{main} to compute the local density of the states $\psi_i^t, i=1,\cdots,N$ in the large $N$ limit. 
The limiting mean local density of state $|\psi_i^t\rangle$ near the energy level $\alpha$ is 
\begin{equation}\label{ldos}
\nu_i(\alpha) = \rho_0(\alpha) \frac{t}{(\alpha-\lambda_i -t H_{\rho_t}(\lambda_i))^2 + t^2 \pi^2 \rho_t(\lambda_i)^2} \,.
\end{equation}  
This last formula \eqref{ldos} is necessarily a probability density function of $\alpha \in \R$, although it is not trivial to check 
that its integral over $\alpha$ is indeed $1$.  
  
\section{Proof of Theorem \ref{main}}  \label{proof}

Following the idea of \cite{sandrine}, we introduce the following quantity, defined for $z\in \C\setminus \R$, as    
\begin{align}
\Theta^g_N(z) &= \frac{1}{N} \sum_{i=1}^N \frac{1}{\lambda_i^t-z} \sum_{j=1}^N  \E[ \langle \psi_i^t|\phi_j\rangle^2]  \,  g(a_j) \notag \\
&= \frac{1}{N} {\rm Tr}\left((M_t-z I)^{-1} g(A) \right)\label{Theta_N.def}
\end{align}
where $g$ is a real valued bounded function on $\R$.  
By convention, $g(A)$ is the diagonal matrix ${\rm Diag}(g(a_1), g(a_2), \cdots, g(a_N))$. 
  
The interesting feature of $\Theta^g_N(z) $ is that, by the Stieltjes inversion formula, we have
\begin{align*}
\Phi_N(\lambda, \alpha) =  \lim_{\eta\rightarrow 0_+} \frac{1}{\pi} \int_{-\infty}^\lambda {\rm Im}\left[\Theta^g_N(\xi+i \eta)\right] d\xi 
\end{align*}
for the particular choice $g(x)= {\bf 1}_{\{x\leq \alpha\}}$.

Thus, the problem is reduced to the study of the convergence of $\Theta^g_N(z)$ when $N\rightarrow \infty$.  
It is plain to deduce Theorem \ref{main} from the following lemma.    
\begin{lemma}\label{lemme.Theta}
Let $g$ be a real valued bounded and continuous function on $\R$. 
Then, as $N\rightarrow +\infty$, we have the following convergence 
\begin{align*}
\Theta^g_N(z) \longrightarrow \Theta^g(z) = \int_0^1 \frac{g(a(x))}{a(x)-z-tG_{\mu_t}(z)} dx
\end{align*}
where $G_{\mu_t}(z)$ is the Stieltjes transform of the limiting spectral distribution $\mu_t$ of the matrix $M_t$.  
\end{lemma}  
  
{\it Proof of Lemma \ref{lemme.Theta}.}  

Using equation \eqref{Theta_N.def} and the definition of the matrix $R_t^N(z)=(M_t-zI)^{-1}$,  
it is straightforward to check that 
\begin{align*}
\Theta^g_N(z) = \frac{1}{N} \sum_{i=1}^N g(a(\frac{i}{N+1})) \, \E[R_t^N(z)_{ii}] \,. 
\end{align*}  
Now, using Theorem \ref{Shlyakhtenko.th} of Shlyakhtenko (see \cite{shlyakhtenko}, and also \cite{pastur,zee}), we know that the complex-valued measure 
$\sigma_t^N$ defined in \eqref{def.sigma} converges weakly to $G_t(z,x) \, dx$.  
Therefore, as $N\rightarrow \infty$, 
\begin{align*}
\Theta^g_N(z) \longrightarrow \int_0^1 g(a(x)) G_t(z,x) dx = \int_0^1 \frac{g(a(x))}{a(x)-z-tG_{\mu_t}(z)} dx\,, 
\end{align*}
using the fixed point equation \eqref{eq.pt-fixe} satisfied by $G_t(z,x)$. The lemma is proved.
 \qed

\section{Initial Semi-circular distribution} \label{semi-circular-initial}
We now turn to the analysis of the particular case where the initial spectral density is a semi-circular distribution. 
 This example is natural in the context of random matrix and free probability theories and we will see that 
 everything can be computed explicitly. 
  
We suppose in this section that the limiting spectral density of the matrix $A$ has a semi-circular shape 
\begin{align}\label{wigner-semi-circle}
\rho_0(x) = \frac{1}{2\pi} \sqrt{4- x^2}\,. 
\end{align}
In order to compute the limiting spectral measure $\mu_t(dx):=\rho_t(x)\, dx$ of the matrix $M(t)$, 
one can write an evolution equation for the Stieltjes transform $G_{\mu_t}(\cdot):\H\to \H$ of $\mu_t$ 
thanks to It\^o's formula (as in Lemma 4.3.12 in \cite{agz}). 
When $N\to \infty$, this {\it Burgers} evolution equation may be written
(see \cite[Proposition 4.3.10]{agz} and \cite[Lemma 3.3.9]{hiai} for a discussion) as
\begin{align}\label{burgers}
\partial_t G_{\mu_t}(z)= G_{\mu_t}(z) \,  \partial _z G_{\mu_t}(z)\,,
\end{align}
and the initial condition is $G_{\mu_0}(z)= (-z+\sqrt{z^2-4})/2$ where the branch of the square root on 
$\C\setminus \R_+$ is such that $\sqrt{-1}=i$. The solution of \eqref{burgers} is easily found
by rewriting \eqref{burgers} in terms of the functional inverse $ B_{\mu_t}$ of the Stieltjes transform $G_{\mu_t}$.   
We easily obtain $B_{\mu_t}(z)= -(t+1) z -1/z$, and for any $z\in \H$, 
\begin{align}\label{Stieltjes-semi-circular-t}
G_{\mu_t}(z) = \frac{-z+\sqrt{z^2-4(1+t)}}{2(1+t)} \,. 
\end{align}
One can recover the limiting spectral  density $\rho_t$ thanks to the Stieltjes inversion formula recalled in \eqref{inversion-formula}, 
\begin{align}\label{semi-circular-t}
\rho_t(x) = \frac{1}{2\pi (1+t)} \sqrt{4 (1+t)- x^2}\,. 
\end{align}
Another method to derive \eqref{semi-circular-t} using the theory of free probability can be found in \cite[ Example 5.3.26]{agz}.

Using \eqref{Stieltjes-semi-circular-t}, we can apply Theorem \ref{main} which 
gives the following expression for the asymptotic overlaps between the $i$-th eigenvector of $M(t)$ and the $j$ th of $A$
with $i\neq j$, 
\begin{align}\label{goe.case}
N \, \E[\langle\psi_i^t|\phi_j\rangle^2] =  
\frac{t}{(a_j-\lambda_i(t))^2  + \frac{t}{1+t} \lambda_i(t) (a_j-\lambda_i(t)) + \frac{t^2}{1+t}}  + o(1)\,. 
\end{align}
Note that, for all $j\neq i$, we have
\begin{equation*}
\lim_{t\to +\infty}\lim_{N\to+\infty} N  \E[\langle \psi_i^t | \phi_j \rangle^2] \longrightarrow 1 \,, 
\end{equation*} 
which means that, after a very long time, 
the vector $|\psi_i(t) \rangle$ has uniform overlaps with the initial eigenvectors $|\phi_j\rangle$ (as one should have expected).
Conversely, the information about a given initial state $|\phi_j\rangle$ is lost when $t\to +\infty$: 
one can not even identify the main components of $|\phi_j\rangle$ in the orthonormal basis of the $|\psi_i(t)\rangle$.

\begin{figure}[h!btp] 
     \center
     \includegraphics[scale=0.75]{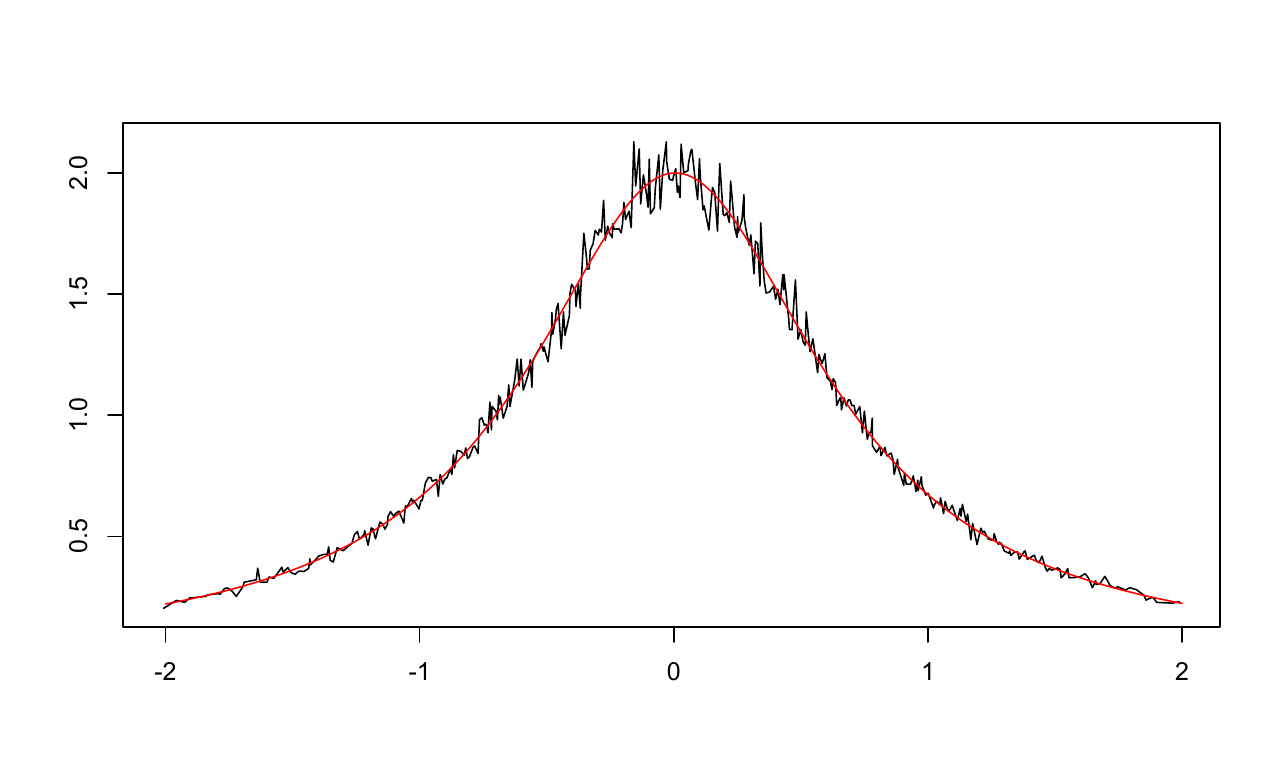}
     \caption{ The black curve is a plot of $N \, \E[\langle \psi_i^t|\phi_j \rangle^2]$, computed empirically with $1000$ samples, 
     as a function of the eigenvalues $a_j$ corresponding 
     to $|\phi_j\rangle$, for $N=400$, $t=1$ and $i=200$. The $200$ th eigenvalue of $A$ is approximately equal to $0$, 
     it is therefore natural to observe 
     the highest value of this curve at this point. The red curve is the theoretical prediction displayed in Eq. \eqref{goe.case}.}\label{fig1}
\end{figure}

\begin{figure}[h!btp] 
     \center
     \includegraphics[scale=0.75]{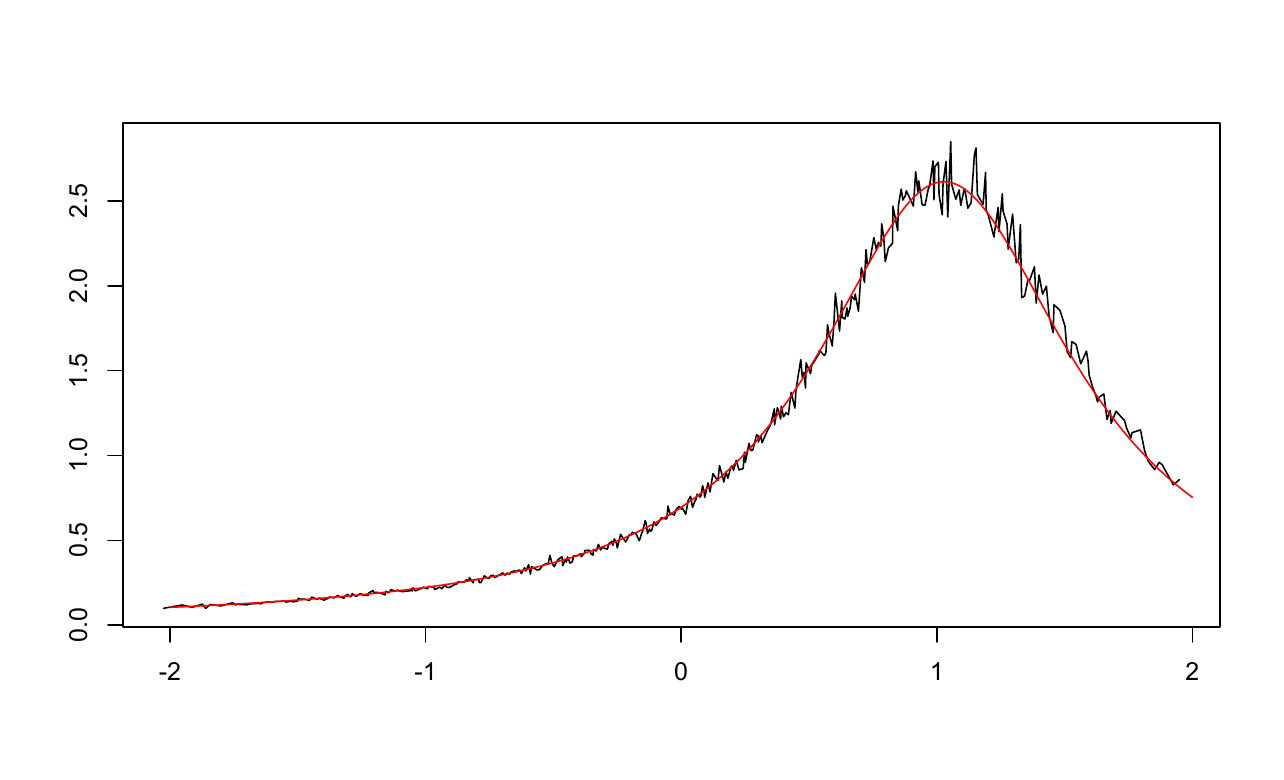}
     \caption{ The black curve is a plot of $N \, \E[\langle \psi_i^t|\phi_j \rangle^2]$, computed empirically with $1000$ samples, 
     as a function of the eigenvalues $a_j$ corresponding 
     to $|\phi_j\rangle$, for $N=400$, $t=1$ and $i=320$. The $320$ th eigenvalue of $A$ is approximately equal to $0.983$, 
     it is therefore natural to observe 
     the highest value of this curve near this point. The red curve is the theoretical prediction displayed in Eq. \eqref{goe.case}.}\label{fig2}
\end{figure}

We have checked the asymptotic expressions \eqref{goe.case} using 
numerical simulations, in the present context where the spectral 
density at the initial time is a semi-circular distribution with variance $1$. 
For this we first sample an initial random symmetric matrix $A$ in the Gaussian Orthogonal Ensemble (GOE) 
(with i.i.d. Gaussian entries up to symmetry with variances $2$ (resp. $1$) on (resp. off) the diagonal),  
of size $N\times N$ with $N= 400$. From Wigner's semi-circular law, we know that
the empirical measure of the eigenvalues of $A$ is close to the semi-circular density 
$\rho_0(x)$ given in \eqref{wigner-semi-circle}.  
For this realization of $A$, we construct a sample of $1000$ matrices $M(t)$ for $t=1$ by adding 
random GOE matrices $H(1)$ to $A$. For all the realizations of $M(1)$ of our sample, we compute the $i$-th
eigenvector with $i=200$ for Fig. \ref{fig1} and $i=320$ for Fig. \ref{fig2} 
and its overlaps with the orthonormal basis of the $|\phi_j\rangle, j=1,\cdots,N$.
We estimate the mean of the overlaps by computing the empirical mean of our sample. 
We finally plot $N \, \E[\langle \psi_i^t|\phi_j \rangle^2]$ as a function of the eigenvalue $a_j$ associated to $|\phi_j\rangle$. 
The agreement between the theoretical and simulated curves is excellent: see Fig. \ref{fig1} and \ref{fig2}.

For Figure \ref{fig1}, we have considered the eigenvector $|\psi_{N/2}(1)\rangle$ 
associated to the median eigenvalue $\lambda_{N/2}(1)$. By symmetry, this eigenvalue is most likely to be near $0$. 
In this particular case where $\lambda_{N/2}(t) \to 0$, Eq. (\ref{goe.case}) can be simplified to read
\begin{align}\label{cauchy}
N \, \E[\langle\psi_{N/2}^t|\phi_j\rangle^2] \approx  \frac{t}{a_j^2  +  t^2}\,, 
\end{align}
which describes a ``Cauchy flight'' in the eigenvalue space of the initial matrix $A$. 
The eigenvector associated to the median eigenvalue $\lambda_{N/2}(t)=0$ 
overlaps the orthonormal basis $|\phi_j\rangle,j=1,\dots,N$, according to a Cauchy distribution in the $A$-eigenvalues space.    
This makes more precise a statement made in \cite{wilkinson,wilkinson2,wilkinson3,vectors} 
in the context of an extreme non adiabatic evolution of a quantum system: the energy is not diffusive but rather performs
a Cauchy Flight.

In fact, if the evolution of the system is such that the elements of the random random Gaussian matrix $M(t)$ 
have a fixed variance, the $i$-th eigenvalue of $M(t)$ is expected to be time independent in the 
large $N$ limit, such that $\lambda_i(t) \approx a_i$. In this case, Eq. (\ref{cauchy}) corresponds (up to 
simple modifications) to Eq. (4.11) of \cite{wilkinson}, with the correspondence $\Delta E= a_j - a_i$. 
However, the correspondence for longer ``times'' $t$ must take into account that with our normalization, 
the semi-circle spectrum itself broadens with time, as given by Eq. \eqref{semi-circular-t}.

\section{Eigenspace stability}\label{stability}

In \cite{vectors}, we investigated the stability of eigenspaces associated to a GOE matrix $A$
when a small GOE matrix $H(t)$ is added. Let us briefly recall the context and main notations of \cite{vectors}.  

Our idea was to study, in the large $N$ limit, the stability of a whole subspace $V_0$ (instead of a single eigenvector as above) 
spanned by a set of consecutive initial eigenvectors $|\phi_k\rangle$ associated to eigenvalues $a_k$ contained 
in a certain interval $[\gamma_-;\gamma_+]$ of 
the Wigner semicircle support $[-2;2]$. 
We then asked the following question: How should one choose a ``larger'' subspace $V_1^t$ spanned by a subset 
 of eigenvectors $|\psi_k^t\rangle$ at time $t$ which would contain the initial subspace $V_0$ up to a small error ?
 To answer this question, we introduced a margin of width $\delta$ and the subspace $V_1^t$ generated by the set of eigenvectors
 $|\psi_k^t\rangle$ associated to eigenvalues $\lambda_k^t$ lying in the interval $[\gamma_- - \delta;\gamma_+ +\delta]$. 
 We then considered the rectangular matrix of overlaps $G_t$ with entries
 \begin{align*}
G_t({ij}) := \langle \psi_i^t |\phi_j \rangle\,. 
\end{align*}
This overlaps matrix is a natural generalization of the scalar product when at least one of 
the dimensions of the two subspaces $V_0$ and $V_1^t$ is greater than one. 
In our setting, the matrix $G_t$ has dimensions $Q \times P$ with 
\begin{equation*}
P = N \int_{\gamma_- }^{\gamma_+} \rho_0, \qquad Q =N \int_{\gamma_- - \delta}^{\gamma_+ + \delta} \rho_t\,,
\end{equation*}
where $\rho_0$ is the Wigner semicircle eigenvalues density of the initial matrix $A$. The labels $i$ and $j$ and the vectors $|\psi_i^t \rangle$ and 
$|\phi_j\rangle$ are respectively indexed by the eigenvalues (in increasing order) $\lambda_i^t$ and $a_j$. 

The $P$ nonzero singular values $1\geq s_1 \geq s_2 \geq \cdots \geq s_P \geq 0$ of the matrix $G_t$ contain 
a meaningful information about the overlap between the two spaces $V_0$ and $V_1^t$. 
For example, the largest singular value $s_1$ indicates that there is a certain linear 
combination of the $Q$ eigenvectors at time $t$ that has a scalar product $s_1$ with a 
certain linear combination of the $P$ initial eigenvectors. If $s_P = 1$, then the initial subspace is 
entirely spanned by the perturbed subspace. If on the contrary $s_1\ll 1$, then the initial 
and perturbed eigenspaces are nearly orthogonal to one another since even the 
largest possible overlap between any linear combination of the original and perturbed eigenvectors is very small. 

A natural way to measure the angle between the subspaces $V_0$ and $V_1^t$ is to compute the quantity 
$v(t)=({\mbox{det}} (G_t^\dagger G_t) )^{1/2}$, which measures the volume of the $P$-dimensional
 parallelepiped spanned by the projection 
of the orthonormal vectors $\phi_1,\cdots,\phi_P$ onto the subspace $V_1^t$. 
In order to get a non trivial limiting $P$-dimensional volume as $P\to \infty$, we need to take a further 
exponent looking at $v(t)^{1/P}$ instead of $v(t)$. 
An even more convenient statistic to measure the angle between the two subspaces $V_0$ and $V_1^t$ 
is in fact obtained by taking the logarithm, $D(V_0,V_1^t) = - \ln ({\mbox{det}} (G_t^\dagger G_t) )^{1/2P}$, 
which can be rewritten as the average of the logarithm of the singular values:
\begin{align*}
D(V_0,V_1^t) = - \frac{1}{P} \sum_{k=1}^P \ln(s_k)\,. 
\end{align*}
This overlap distance $D$ and the overlap matrix $G_t$ already 
appeared in the literature on the ``Anderson orthogonality catastrophe" (see e.g. \cite{anderson,ullmo}).

Using perturbation theory, we showed in \cite{vectors} that this overlap distance $D(V_0,V_1^t)$, behaves when 
$N \to +\infty$ with $t:= t_N$ such that $Nt_N\to 0$, as:
\begin{align}\label{eq.D}
\E[D(V_0,V_1^t)] = \frac{t}{2\int_{\gamma_-}^{\gamma_+} \rho_0} \int_{\gamma_-}^{\gamma_+} dx  \int_{y\not \in [\gamma_--\delta,\gamma_++\delta]} dy \, \frac{\rho_0(x) \rho_0(y)}{(x-y)^2} + o(t)\,. 
\end{align}
The fixed parameter $\delta>0$ permits us to truncate the singularity induced by 
 pseudo collisions at the edge of the intervals $[\gamma_-,\gamma_+]$. 
Supported by 
convincing numerical evidence \cite{vectors}, 
this formula \ref{eq.D} was conjectured to hold true in the semi-perturbative regime $N\to +\infty$ 
with $t_N \to 0$ but not necessarily $N t_N \to 0$. Unfortunately we were unable at the time 
to find analytical arguments to sustain our claim in this semi-perturbative regime. 

Our new results obtained in this paper provide us good tools to fill in this gap and prove that Eq. \eqref{eq.D} is 
indeed correct in the semi-perturbative regime $t:=t_N\to 0, N\rightarrow +\infty$.  

Let us first remark that $D(V_0,V_1^t)= - \ln(\det(G_t^\dagger G_t))/(2P)$ where $G_t^\dagger$ is the Hermitian conjugate of $G_t$.   
We start by computing the entries of the matrix $G_t^\dagger G_t$. 
Using the fact that the $|\psi_k^t\rangle, k=1,\cdots,N$ form an orthonormal family of $\R^N$, we have 
for all $a_i \in [\gamma_-; \gamma_+]$,
\begin{align*}
(G_t^\dagger G_t)_{ii} = \sum_{\lambda_k^t \in [\gamma_- - \delta; \gamma_+ +\delta]} \langle \psi_k^t|\phi_i\rangle^2
= 1 -  \sum_{\lambda_k^t \not\in [\gamma_- - \delta; \gamma_+ +\delta] } \langle \psi_k^t|\phi_i\rangle^2\,.
\end{align*}
Using \eqref{perturb.2.n}, we see that, as $t_N\to 0$ and $N\to \infty$,
\begin{align}
\E[(G_t^\dagger G_t)_{ii}] & = 1 -  \frac{t}{N} \sum_{a_k \not \in [\gamma_- - \delta; \gamma_+ +\delta] } \frac{1}{(a_i-a_k)^2} + o(\frac{1}{N}) +o(t) \notag \\
&\sim_{N \to \infty} 1 -  t \int_{y\not \in [\gamma_- -\delta,\gamma_+ +\delta]} dy\, \frac{\rho_0(y)}{(a_i - y)^2} + o(t)\,.  \label{diag.subspace}
\end{align} 
The non diagonal elements, i.e. indexed by $a_i \neq a_j \in [\gamma_-;\gamma_+]$, can also be computed as
\begin{align*}
(G_t^\dagger G_t)_{ij} =  \sum_{\lambda_k^t \in [\gamma_- - \delta; \gamma_+ +\delta]} \langle \psi_k^t|\phi_i\rangle  \langle \psi_k^t|\phi_j\rangle 
= - \sum_{\lambda_k^t \not \in [\gamma_- - \delta; \gamma_+ +\delta]} \langle \psi_k^t|\phi_i\rangle  \langle \psi_k^t|\phi_j\rangle
\end{align*}
where, in the second line, we have used the orthogonality of $|\phi_i\rangle$ and $|\phi_j\rangle$ which implies that 
$\sum_k \langle \psi_k^t|\phi_i\rangle  \langle \psi_k^t|\phi_j\rangle = 0$. The expectations of those terms can thus be estimated 
as  $t:=t_N\to 0$ and 
 $N\to \infty$, via the Cauchy-Schwarz inequality, as 
\begin{align}
\E[(G_t^\dagger G_t)_{ij}] &\leq  \sum_{a_k \not \in [\gamma_- - \delta; \gamma_+ +\delta] }   \E[\langle \psi_k^t|\phi_i\rangle]^{1/2}
\E[\langle \psi_k^t|\phi_j\rangle]^{1/2} \notag \\ 
& \sim_{N\to +\infty} \frac{t}{N} \sum_{a_k \not \in [\gamma_- - \delta; \gamma_+ +\delta] }  \frac{1}{(a_i-a_k)(a_j-a_k)} + o(\frac{1}{N}) + o(t) \\
&\sim_{N \to \infty}  t   \int_{y\not \in [\gamma_- -\delta,\gamma_+ +\delta]} dy \, \frac{\rho_0(y)}{(a_i - y)(a_j-y)} +o(t) \label{non.diag.subspace}  \,. 
\end{align}
Note that in both equations \eqref{diag.subspace} and \eqref{non.diag.subspace}, we have $a_i, a_j \in [\gamma_-; \gamma_+]$
and $a_k \not\in [\gamma_- - \delta;\gamma_+ + \delta]$,
so that $a_i$ (resp. $a_j$) and $a_k$ remain at macroscopic distance $> \delta$ and formula \eqref{perturb.2.n} applies. 
Besides the integrals in \eqref{diag.subspace} and \eqref{non.diag.subspace} 
are perfectly well defined due to the introduction of the margin $\delta>0$. 

Thus, if $t:=t_N\to 1$ and $N\to \infty$, the determinant of $G_t^\dagger G_t$ 
can be approximated to leading order in $t_N$ as the product of the diagonal terms 
(the other contribution are negligible compared to $t_N$). 
We thus have, doing a further linearization when $t_N\to 0$, 
\begin{align*}
- \frac{1}{2P} \E[\ln( \det(G_t^\dagger G_t))] &= \frac{t}{ P } \sum_{a_i \in [\gamma_-;\gamma_+]}  
\int_{y\not \in [\gamma_- -\delta,\gamma_+ +\delta]} dy \, \frac{\rho_0(y)}{(a_i - y)^2} + o(\frac{1}{N})+ o(t) \\
& \sim_{N\to +\infty} \frac{t}{\int_{\gamma_-}^{\gamma_+} \rho_0} \int_{\gamma_-}^{\gamma_+} dx \int_{y\not \in [\gamma_--\delta,\gamma_++\delta]} dy \, \frac{\rho_0(x) \rho_0(y)}{(x- y)^2} + o(t) \,.
\end{align*}
This is our proof that \eqref{eq.D} is valid in the second semi-perturbative regime. 

The reader may wonder how to extend formula \eqref{eq.D} in the non perturbative regime, i.e. for arbitrary values of $t$. 
This question is clearly more difficult as one would need to understand 
the convergence of the non diagonal terms of the matrix $G_t^\dagger G_t$ 
in the large $N$ limit, which are no longer negligible in the determinant expansion.

\end{document}